\documentclass[12pt,a4paper]{amsart}

\usepackage{graphicx}
\usepackage{amsfonts,amsmath,latexsym,amssymb,amsthm}

\newcounter{theoremcounter}
\newcounter{lemmacounter}

\newcounter{dummycounter}

\newcounter{conjcounter}

\newcounter{emptycounter}
\newcounter{defcounter}

\newtheorem{theorem}[theoremcounter]{Theorem}

\newtheorem{lemma}[lemmacounter]{Lemma}
\newtheorem{conjecture}[conjcounter]{Conjecture}

\newtheorem{definition}[defcounter]{Definition}

\numberwithin{equation}{section}
\numberwithin{lemmacounter}{section}
\numberwithin{propcounter}{section}
\numberwithin{corcounter}{section}
\numberwithin{conjcounter}{section}
\numberwithin{theoremcounter}{section}
\numberwithin{probcounter}{section}

\newcounter{eqncounter}

\numberwithin{equation}{eqncounter}

\def\X{T}
\def\Mcp{\mathcal{M}_{K}^{(cp)}}
\def\MFcp{\mathcal{M}_{\L}^{(cp)}}
\def\Mred{\mathcal{M}_{K}^{(red)}}
\def\Mncp{\mathcal{M}_{K}^{(ncp)}}
\def\Mall{\mathcal{M}_{K}}
\def\MFall{\mathcal{M}_{\L}}
\def\grmn{\mathcal{M}_K(n,\X)}
\def\grmnm1{\mathcal{M}_K(n-1,\X)}
\def\grmnnew{\mathcal{M}_K(n,X^n)}
\def\grmnmnew1{\mathcal{M}_K(n-1,X^n)}

\newcommand\pw{\mathfrak{p}}
\newcommand\f{l}

\newcommand\coll{\mathcal{C}}

\newcommand\ALS{ALS}

\newcommand\e{e}

\newcommand\M{M}
\newcommand\enMK{\en'_K}
\newcommand\enML{\en'_F}
\newcommand\enM{\en'}

\newcommand\IR{\mathbb R}
\newcommand\IC{\mathbb C}
\newcommand\IZ{\mathbb Z}

\newcommand\IP{\mathbb P}
\newcommand\IQ{\mathbb Q}

\newcommand\en{\mathcal{N}}
\newcommand\hen{H_\mathcal{N}}
\newcommand\henMK{H_{\mathcal{N}_K'}}
\newcommand\henML{H_{\mathcal{N}_F'}}
\newcommand\henM{H_\mathcal{N'}}
\newcommand\balf{{\mbox{\boldmath $\alpha$}}}
\newcommand\vnull{{\mbox{\boldmath $0$}}}
\newcommand\Hom{\text{\rm{Hom}}}

\newenvironment{rproof}{\addvspace{\medskipamount}\par\noindent{\it Proof.\/}}
{\unskip\nobreak\hfill$\Box$\par\addvspace{\medskipamount}}
\newcommand\ord{\mathop{\rm ord}\nolimits}

\newcommand\Oseen{{\mathcal{O}}}

\renewcommand{\L}{{F}}

\newcommand\Qbar{\overline{\IQ}}

\newcommand\vx{{\bf x}}
\newcommand\vy{{\bf y}}
\newcommand\vz{{\bf z}}

\newcommand\henK{H_{{\mathcal{N}_{K}}}}

\newcommand\Ce{D}
\newcommand\Da{D}

\begin{document}\baselineskip=17pt
\title{On number fields with nontrivial subfields}

\author{Martin Widmer}

\address{Institut f\"ur Mathematik\\ 
Technische Universit\"at Graz\\ 
Steyrergasse 30/II\\
A-8010 Graz\\ 
Austria}

\email{widmer@tugraz.at}

\date{October 5, 2009}



\begin{abstract}
What is the probability for a number field of composite degree $d$ to have a nontrivial subfield? As the reader might expect the answer heavily depends 
on the interpretation of probability. We show that if the fields are enumerated by the smallest height of their generators 
the probability is zero, at least if $d>6$. This is in contrast to what one expects when the fields are enumerated by the discriminant.
The main result of this article is an estimate for the number of algebraic numbers of degree $d=\e n$ and bounded height 
which generate a field that contains an unspecified subfield of degree $\e$. If $n>\max\{\e^2+\e,10\}$ we get the correct asymptotics as the height tends to infinity.
\end{abstract}

\maketitle

\section{Introduction and results}

The most natural way to enumerate number fields of fixed degree is probably by their discriminant $\Delta$ or the absolute value thereof.
For a positive integer $d$ let $\Delta(d,X)$ be the number of field extensions $\L$ of $\IQ$ of degree $d$
in an algebraic closure $\Qbar$ with $|\Delta_\L|\leq X$.
The asymptotics are predicted by a classical
conjecture, possibly due to Linnik (see e.g. \cite{56}), but proved only for 
degree $d=2,3,4,5$.
\begin{conjecture}
Suppose $d>1$.
Then there exists a positive constant $c_{d}$ such that as $X$ tends to infinity
\begin{alignat*}1
\Delta(d,X)=c_{d}X+o(X).
\end{alignat*}
\end{conjecture}
Linnik's Conjecture is usually stated in a more general form which asserts that for any number field $K$ the number of
field extensions $\L$ of $K$ of relative degree $n$ satisfying $|\Delta_\L|\leq X$ is given by
$c_{K,n}X+o(X)$ for a positive constant $c_{K,n}$.\\

Let $G$ be a subgroup of the symmetric group $S_{d}$
containing a subgroup of index $d$. Malle \cite{57} has given
conjectural asymptotics for $\Delta_G(d,X)$, the number of fields in $\Qbar$
of degree $d$
whose Galois closure has Galois group isomorphic to $G$
and whose absolute value of the discriminant is not larger than $X$. 
Kl\"uners \cite{Kluners2005} found counterexamples to Malle's conjecture but a slight adjustment of the conjecture proposed
by T\"urkelli \cite{Turkelli2008} seems promising. But once
again this is proved only in very special cases. Bhargava's work
\cite{57}
implies $\Delta_{S_4}(4,X)\thicksim \lambda X$ for 
\begin{alignat*}1
\lambda=\frac{5}{6}\prod_p\left(1+\frac{1}{p^2}-\frac{1}{p^3}-\frac{1}{p^4} \right)=1.01389....
\end{alignat*}
And according to Cohen, Diaz y Diaz and Olivier \cite{59}
the number with Dihedral group $\Delta_{D_4}(4,X)$ is $\thicksim \mu X$
where $\mu=0.1046520224...$.
A quartic field has a quadratic subfield if and only if
its Galois closure is $D_4$ or an abelian group of order four.
Bailey \cite{Bailey1980} and Wong \cite{Wong1999} have shown that 
$\Delta_G(4,X)=o(X)$ for $G=A_4$ and abelian groups $G$ of order four. Thus when we enumerate the quartic fields by the absolute value of their discriminant
the probability that a quartic field has a quadratic subfield
is the positive number
\begin{alignat*}1
\frac{\mu}{\mu+\lambda}=0.09356....
\end{alignat*}
Suppose the (generalized) Linnik Conjecture is true. We fix a number field $K$ of degree $\e$ and then we count extensions $\L$ of $K$
of relative degree $n$ satisfying $|\Delta_\L|\leq X$. In this way we conclude that the lower density for the set of fields of degree $d=\e n$ that contain a subfield of degree $\e$ is positive; of course
here density is understood with respect to the absolute value of the discriminant. Hence when enumerated by the absolute value of the discriminant the (``lower'') probability that a field of degree $\e n$
has a subfield of degree $\e$ remains positive, subject to the (generalized) Linnik Conjecture. \\

This is in stark contrast to
the situation when one enumerates by the following, also classical, invariant
\begin{alignat*}1
\pi(\L)=\inf_{\alpha \atop \IQ(\alpha)=\L}|D_{\alpha}|.
\end{alignat*}
Here $D_{\alpha}$ is the unique minimal polynomial of $\alpha$
in $\IZ[x]$ with positive leading coefficient and coprime coefficients and $|D_{\alpha}|$ denotes the maximum norm
of the coefficient vector. The quantity $|D_{\alpha}|$ is sometimes referred to as the naive height of $\alpha$. We define the counting function $\pi(\e, n,X)$ as the number of fields $\L\subseteq \Qbar$ of degree $\e n$
that contain a subfield of degree $\e$ and satisfy $\pi(\L)\leq X$.\\

In this note we shed some light on the distribution of number fields by 
counting generators. Let $H$ be the absolute multiplicative Weil height (or briefly the height) on $\Qbar$, as defined in \cite[p.16]{BG}.
A result of Masser and Vaaler (\cite[Theorem]{37}) gives the asymptotics for the number
of generators of degree $\e n$ with bounded height. We extend Masser and Vaaler's result by estimating
$Z(\e,n,X)$ which counts the numbers with height at most $X$
generating a field of degree $\e n$ that contains a subfield of  degree $\e$
\begin{alignat*}1
Z(\e,n,X)=|\{\alpha\in\Qbar;[\IQ(\alpha):\IQ]=\e n,
\IQ(\alpha)\text{ contains a field of degree $\e$, } H(\alpha)\leq X \}|.
\end{alignat*}
Our first result is a simple by-product of the proof of our main result Theorem \ref{th1} combined with a result of Schmidt, and gives an upper bound for 
$Z(\e,n,X)$.
\begin{theorem}\label{th4}
With $c=n\cdot 2^{\e(n^2+n\e+2\e+n+13)+n^2+10n}$ and $X>0$ we have
\begin{alignat*}1
Z(\e,n,X)\leq cX^{\e n(n+\e)}.
\end{alignat*}
\end{theorem}
The invariant $\delta(\L)=\inf\{H(\alpha); \L=\IQ(\alpha)\}$ plays a crucial role in the proofs. 
If $\alpha$ is an algebraic number of degree $\e n$ then $H(\alpha)^{\e n}=M(D_\alpha)$ where $M$ denotes the Mahler measure (see \cite[p.22]{BG} or \cite[p.434]{1} for a definition). A crude estimate comparing $M(D_\alpha)$ and $|D_\alpha|$ gives
\begin{alignat}1\label{heightmaxnorm}
(2^{-1}H(\alpha))^{\e n}\leq|D_{\alpha}|\leq (2H(\alpha))^{\e n}
\end{alignat}
and hence
\begin{alignat*}1
(2^{-1}\delta(\L))^{\e n} \leq \pi(\L) \leq (2\delta(\L))^{\e n}.
\end{alignat*}
We therefore conclude from Theorem \ref{th4}
\begin{alignat*}1
\pi(\e,n,X)\leq c\cdot2^{\e n(n+\e)}X^{n+\e}.
\end{alignat*}
On the other hand Corollary 5.1 in \cite{art2} yields
\begin{alignat*}1
\pi(1,\e n,X)\geq C_{\e n}X^{\e n-1}
\end{alignat*}
for a positive constant $C_{\e n}$ and $X\geq X_0(\e n)$.
Combining these two estimates we find:
when ordered by the invariant $\pi$ the 
probability that a field $\L$ of degree $\e n$ has a subfield 
different from $\IQ$ and $\L$ is zero, at least for $\e n>6$.\\

Another consequence of Theorem \ref{th4}
concerns polynomials with certain Galois groups.
Let $f$ in $\IZ[x]$ be irreducible of degree $\e n$. 
Since Van der Waerden \cite{vanderWaerden1934} it is known that almost all
polynomials $f$ have the full symmetric group $S_{\e n}$
as Galois group when enumerated by the maximum norm of the coefficient vector.
That is any root $\alpha$ of $f$ generates a field $\L=\IQ(\alpha)$
whose Galois closure $\L_G$ has Galois group $S_{\e n}$ over
$\IQ$. The group corresponding to $\L$ is some $S_{\e n-1}$.
It is easy to see that there is no group lying strictly
between these two groups. This means that $\L/\IQ$ has no proper intermediate field in this case. Van der Waerden's result can be further quantified through sharpenings of the Hilbert Irreducibility Theorem. A general version due to S.D. Cohen (\cite[Theorem 2.1]{72})
gives an upper bound of order $X^{\e n+1/2}\log X$ for the number 
of exceptional polynomials. Gallagher and Dietmann \cite{Dietmann2006} improved the exponent $\e n+1/2$ for $\e n=4$. It is likely that the exponent $\e n+1/2$ can always be improved but this might be hard to achieve in general. 
However, under the stronger condition
that there exists a proper intermediate field Theorem \ref{th4} in combination with (\ref{heightmaxnorm}) tells us that the exponent $\e n +1/2$ can be reduced to $\e n/2+2$.\\

So much for the consequences of the proof of our main result. We now come to the main result itself.
As already mentioned it asymptotically estimates the counting function $Z(\e,n,X)$ as the height bound $X$ tends to infinity.
To state the result we have to introduce further notation.
In \cite{1} Masser and Vaaler defined the following two quantities 
\begin{alignat*}1
V_{\IR}(n)=(n+1)^l\prod_{i=1}^{l}\frac{(2i)^{n-2i}}{(2i+1)^{n+1-2i}}
\end{alignat*}
where $l=[(n-1)/2]$ and the empty product is interpreted as $1$ and
\begin{alignat*}1
V_{\IC}(n)=\frac{(n+1)^{n+1}}{((n+1)!)^{2}}.
\end{alignat*}
These formulae give the volumes of the unit balls in $\IR^{n+1}$ and $\IC^{n+1}$ with respect to 
the Mahler measure distance function and have been calculated by Chern and Vaaler in \cite{43}.
We also need the Schanuel constant $S_K(n)$ for a number field $K$, defined as follows
\begin{alignat}1\label{Schanuelconst}
S_K(n)=\frac{h_KR_K}{w_K\zeta_K(n+1)}
\left(\frac{2^{r_K}(2\pi)^{s_K}}{\sqrt{|\Delta_K|}}\right)^{n+1}
(n+1)^{r_K+s_K-1}.
\end{alignat}
Here $h_K$ is the class number, $R_K$ the regulator,
$w_K$ the number of roots of unity in $K$, $\zeta_K$
the Dedekind zeta-function of $K$, $\Delta_K$ the discriminant,
$r_K$ is the number of real embeddings of $K$ and $s_K$ is
the number of pairs of distinct complex conjugate embeddings of $K$.\\

All fields are considered to lie in a fixed algebraic closure $\Qbar.$
It will be convenient to use Landau's $O$-notation.
For non-negative real functions $f(X), g(X), h(X)$ we say that
$f(X)=g(X)+O(h(X))$ as $X>X_0$ tends to infinity if there is a constant $C_0$ such that
$|f(X)-g(X)|\leq C_0h(X)$ for each $X>X_0$.
Now we can state the main result.
\begin{theorem}\label{th1}
Suppose $n> \max\{\e^2+\e,10\}$.
Then as $X>0$ tends to infinity we have
\begin{alignat}1
\label{eqth1}
Z(\e,n,X)=\left(\sum_K nV_{\IR}(n)^{r_K}V_{\IC}(n)^{s_K}S_K(n)\right)X^{\e n(n+1)}
+O(X^{\e n(n+1)-n}),
\end{alignat}
where the sum runs over all number fields of degree $\e$
and the implied constant in the $O$-term depends only on $\e$ and $n$.
\end{theorem}
The above theorem states implicitly, subject to the constraints on $\e$ and $n$, that the
sum on the right-hand side of (\ref{eqth1}) converges.
Notice that by Masser and Vaaler's Theorem \cite{37} (or its generalization from $\IQ$ to arbitrary ground fields in \cite{1}) 
\begin{alignat*}1
Z(1,\e n,X)=Z(\e n,1,X)= \e nV_{\IR}(\e n)S_{\IQ}(\e n)X^{\e n(\e n+1)}
+O(X^{(\e n)^2}\mathfrak{L})
\end{alignat*}
where $\mathfrak{L}$ is defined in Theorem \ref{th3}.
So for instance the asymptotics for the numbers of degree $22$ involve $X^{506}$ whereas those for the numbers that generate a field which contains a quadratic subfield involve only $X^{264}$.\\

If each divisor $>1$ of $n$ is larger than $\e$
we can relax the constraints on $\e$ and $n$.
\begin{theorem}\label{th2}
Suppose $\f>1$ and $\f|n$ implies $\f>\e$ 
and suppose $n>\max\{6\e-6,10\}$. Then as $X>0$ tends to infinity we have
\begin{alignat*}1
Z(\e,n,X)=\left(\sum_K nV_{\IR}(n)^{r_K}V_{\IC}(n)^{s_K}S_K(n)\right)X^{\e n(n+1)}
+O(X^{\e n(n+1)-n})
\end{alignat*}
where the sum runs over all number fields of degree $\e$. The implied constant in the $O$-term depends only on $\e$ and $n$.
\end{theorem}
Our proof strategy for Theorem \ref{th1} can be roughly (and oversimplified) described as follows. First fix a field $K$ of degree $\e$ and count those
numbers having degree $n$ over $K$ and degree $\e n$ over $\IQ$.
Combining ideas of Masser and Vaaler from \cite{1} and of the author's works \cite{art1} and \cite{art2} this can be achieved 
by counting monic polynomials $x^n+\alpha_1x^{n-1}+\cdots+\alpha_n$ in $K[x]$ with $K=\IQ(\alpha_1,...,\alpha_n)$ and with bounded Mahler measure. For the error term one has to take into account the reducible polynomials and
also the polynomials irreducible over $K$ but reducible over the
Galois closure of $K$. Then we sum these estimates over all fields $K$ of degree $\e$.
This requires that the emerging error terms converge when summed over all fields $K$. 
The error terms are expressed using the invariant $\delta(K)$, because they have better summatory properties than
the discriminant. \\

We can use the same ideas to prove asymptotic results for 
\begin{alignat*}1
Z(\e,m,n,X)=|\{&\alpha\in\Qbar;[\IQ(\alpha):\IQ]=\e m n, H(\alpha)\leq X,\\
&\IQ(\alpha)\text{ contains a field of degree $\e$ and a field of degree $\e m$}\}|.
\end{alignat*}
We state just one particularly simple result.
\begin{theorem}\label{th5}
Suppose $\f>1$ and $\f|n$ implies $\f>\e m$ 
and suppose $n>\max\{6\e m-6,10\}$. Then as $X>0$ tends to infinity we have
\begin{alignat*}1
Z(\e,m,n,X)=\left(\sum_K nV_{\IR}(n)^{r_K}V_{\IC}(n)^{s_K}S_K(n)\right)X^{\e m n(n+1)}
+O(X^{\e mn(n+1)-n})
\end{alignat*}
where the sum runs over all number fields of degree $\e m$ that contain a subfield of degree $\e$.
\end{theorem}
Notice that under the above conditions on $\e,m$ and $n$ the functions $Z(1,\e m,n,X)$ and $Z(\e,m,n,X)$ both have order of magnitude $X^{\e m n(n+1)}$ whereas $Z(1,1,\e mn,X)$ has order of magnitude $X^{\e m n(\e mn+1)}$.\\

Let us mention one final by-product of the proof of Theorem \ref{th1}. We obtain a version of
the Theorem in \cite{1} with a particularly good error term regarding the ground field $K$ 
under the necessary condition that we exclude those numbers
that have also degree $n$ over a proper subfield $k$ of $K$.
\begin{theorem}\label{th3}
Let $K$ be a number field of degree $\e$.
Then as $X>0$ tends to infinity the number of elements $\beta$ in $\Qbar$ with
\begin{alignat}1
\nonumber&[K(\beta):K]=n, \\ 
\label{theorem3cha4cond2}
&k\subseteq K \text{ and }[k(\beta):k]=n \Longrightarrow k=K, \\ 
\nonumber &H(\beta)\leq X
\end{alignat}
is 
\begin{alignat*}1
nV_{\IR}(n)^{r_K}V_{\IC}(n)^{s_K}S_K(n)X^{\e n(n+1)}
+O(\delta(K)^{-\frac{\e}{2}(n-\max\{4\e-8,2\e-3\})+1.1}X^{\e n(n+1)-n}\mathfrak{L})
\end{alignat*}
where $\mathfrak{L}=1$ unless $\e n=1$ or $\e n=2$ in which case
$\mathfrak{L}=\log (X+2)$. The constant in $O$ depends only on $\e$ and $n$.
\end{theorem}
If $\e$ and $n>\max\{4\e-8,2\}$ are fixed then the constant in the error term goes rapidly to zero as the fields $K$ become more complicated.
The additive constant $1.1$ in the exponent on $\delta(K)$ has no particular significance and could be replaced by any other value $>1$.\\
For $\e=1$ or $n=1$ Theorem \ref{th4} is covered by Schmidt's Theorem in \cite{22}.
The cases $\e=1$ in Theorem \ref{th1}, Theorem \ref{th2} and Theorem \ref{th3} are all covered by Masser and Vaaler's Theorem in \cite{37} and the case $n=1$ in Theorem \ref{th3} counts generators $\alpha \in K$ with bounded height and thus is
covered by a special case of Corollary 3.2 in \cite{art1} (which we cite as Theorem \ref{prop2} in Section \ref{introchap4}).
Finally the cases $\e=1$ or $m=1$ in Theorem \ref{th5} are covered by Theorem \ref{th2}.
We emphasize that our work neither gives a proof of Schmidt's nor a new proof of Masser and Vaaler's result but 
rather uses their method and ideas in combination with the work done in \cite{art1} and \cite{art2} to extend these results.\\
Throughout this article $X$ and $T$ denote positive real numbers.

\section*{Acknowledgements}
Although the presented work is not included in my Ph.D. thesis the major part has been 
carried out during my graduate studies.
I am indebted to my Ph.D. adviser David Masser for his generosity in sharing his thoughts and ideas, in particular
for pointing out to me that the results from \cite{art2} might yield results in the style of Theorem \ref{th1}.
I also would like to thank Wolfgang Schmidt for interesting and fruitful discussions.
This work was financially supported by the Swiss National Science Foundation.

\section{Reformulation of Theorem \ref{th1} step one}
Let $K$ be a number field of degree $\e$. We define
\begin{alignat*}1
Z_K(\e,n,X)=|\{\beta\in \Qbar;[\IQ(\beta):\IQ]=\e n, [K(\beta):K]=n,
H(\beta)\leq X\}|.
\end{alignat*}
If $\beta \in \Qbar$
with $[\IQ(\beta):\IQ]=\e n$ and $\IQ(\beta)$
contains the field $K$ of degree $\e$ then $[K(\beta):K]=n$.
Therefore 
\begin{alignat}1
\label{ZmnXub2}
Z(\e,n,X)\leq \sum_K Z_K(\e,n,X),
\end{alignat}
where $K$ runs over all fields of degree $\e$.  On the other hand
if $\beta$ is in $\Qbar$ with $[K(\beta):K]=n$ and 
$[\IQ(\beta):\IQ]=\e n$ then $\IQ(\beta)$ contains
the field $K$ of degree $\e$. However, some elements $\beta$ may be counted for several different fields $K$
on the right-hand side of (\ref{ZmnXub2}).
To keep track of these multiply counted numbers
we have to introduce two further quantities.
\begin{alignat*}1
\overline{Z}(\e,n,X)=&\\
&|\{\beta\in \Qbar;[\IQ(\beta):\IQ]=\e n,\\
&\IQ(\beta)\text{ contains more than one field of degree $\e$},H(\beta)\leq X\}|,\\
 \overline{Z}_K(m,n,X)=&\\
&|\{\beta\in \Qbar;[\IQ(\beta):\IQ]=\e n, [K(\beta):K]=n,\\
&\IQ(\beta)\text{ contains more than one field of degree $\e$},H(\beta)\leq X\}|.
\end{alignat*}
For all $\e,n$ we have 
\begin{alignat}1
\label{3gl2}
Z(\e,n,X)=\sum_K \left(Z_K(\e,n,X)-\overline{Z}_K(\e,n,X)¸\right)
+\overline{Z}(\e,n,X).
\end{alignat}
where $K$ runs over all fields of degree $\e$.
Moreover
\begin{alignat}1
\label{3ugl1}
\overline{Z}(\e,n,X)\leq \sum_K \overline{Z}_K(\e,n,X)\leq
2^{\e n}\overline{Z}(\e,n,X).
\end{alignat}
The first inequality is obvious; the second one holds
because every field of degree $\e n$ contains at most $2^{\e n}$ 
subfields.\\
Now suppose $\IQ(\beta)$ contains more than one subfield
of degree $\e$. So the compositum of two different
subfields of degree $\e$ lies in $\IQ(\beta)$. But this
compositum has degree $\f\e$ where $\f\mid n$ and
$\f \in  \{2,3,...,\e\}$. Hence by (\ref{ZmnXub2})
\begin{alignat*}1
\overline{Z}(\e,n,X)\leq \sum_{\f|n \atop{1<\f\leq \e}}
Z(\f\e,n/\f,X)\leq \sum_{\f|n \atop{1<\f\leq \e}}
\sum_{\L \atop{[\L:\IQ]=\f\e}} Z_\L(\f\e,n/\f,X).
\end{alignat*}
Together with (\ref{3gl2}) and (\ref{3ugl1}) we get
\begin{alignat}1
\label{3gl3}
Z(\e,n,X)=\sum_{K \atop [K:\IQ]=\e} Z_K(\e,n,X)
+O\left(\sum_{\f|n \atop{1<\f\leq \e}}
\sum_{\L \atop{[\L:\IQ]=\f\e}} Z_\L(\f\e,n/\f,X)\right)
\end{alignat}
The sums in (\ref{3gl3}) can essentially be reduced to the counting of projective points $P$ in $\IP^n$
of degree $\e$ with $\hen(P)\leq X$ for a certain adelic-Lipschitz height $\hen$. The next section is devoted to 
the basic definitions of this concept and the necessary results to derive the statements of this article.

\section{Adelic-Lipschitz systems and adelic-Lipschitz heights}
This section is (in fact in a more general form) contained in \cite{art2}.
However, for convenience of the reader we recall the general concept
of an adelic-Lipschitz system and its basic definitions.
\subsection{Adelic-Lipschitz systems on a number field}\label{subsecdefALS}
Let $r$ be the number of  
real embeddings and $s$ the number of pairs of complex conjugate embeddings
of $K$ so that $\e=r+2s$.
Recall that $M_K$ denotes the set of places of $K$.
For every place $v$ we fix a
completion $K_v$ of $K$ at $v$ and we write $d_v=[K_v:\IQ_v]$ with $\IQ_v$ being the completion with respect to the place that extends to $v$.
A place $v$ in $M_K$ corresponds either to a non-zero prime ideal $\pw_v$
in the ring of integers $\Oseen_K$ or to an embedding
$\sigma$ of $K$ into $\IC$.
If $v$ comes from a prime ideal we call $v$
a finite or non-archimedean place and denote this by $v\nmid \infty$ and if $v$ corresponds to
an embedding we say $v$ is an infinite or archimedean
place and denote this by $v\mid \infty$.
For each place in $M_K$ we choose a representative
$|\cdot|_v$,
normalized in the following way:
if $v$ is finite and $\alpha\neq 0$ we set by convention
\begin{alignat*}3
|\alpha|_{v}=N\pw_v^{-\frac{\ord_{\pw_v}(\alpha\Oseen_K)}{d_v}}
\end{alignat*}
where $N\pw_v$ denotes the norm of $\pw_v$
from $K$ to $\IQ$ and $\ord_{\pw_v}(\alpha\Oseen_K)$
is the power of $\pw_v$ in the prime ideal decomposition
of the fractional ideal $\alpha\Oseen_K$.
Moreover we set
\begin{alignat*}3
|0|_{v}=0.
\end{alignat*}
And if $v$ is infinite and corresponds to an embedding $\sigma:K \hookrightarrow \IC$ we define
\begin{alignat*}3
|\alpha|_{v}=|\sigma(\alpha)|.
\end{alignat*}
The value set of $v$, $\Gamma_v:=\{|\alpha|_v;\alpha \in K_v\}$
is equal to $[0,\infty)$ if $v$ is archimedean,
and to
\begin{alignat*}3
\{0,(N\pw_v)^{0},(N\pw_v)^{\pm 1/d_v},(N\pw_v)^{\pm 2/d_v},...\}
\end{alignat*}
if $v$ is non-archimedean.
For $v \mid \infty$ we identify $K_v$ with $\IR$ or
$\IC$ respectively and we identify $\IC$ with
$\IR^2$ via $\xi\longrightarrow (\Re(\xi),\Im(\xi))$
where we used $\Re$ for the real and $\Im$ for the
imaginary part of a complex number.\\
For a vector $\vx$ in $\IR^n$ we write $|\vx|$ for the euclidean length of $\vx$.
\begin{definition}
Let $\M$ and $\Da>1$ be positive integers and let $L$ be a non-negative real.
We say that a set $S$ is in Lip$(\Da,\M,L)$ if 
$S$ is a subset of $\IR^\Da$, and 
if there are $\M$ maps 
$\phi_1,...,\phi_M:[0,1]^{\Da-1}\longrightarrow \IR^\Da$
satisfying a Lipschitz condition
\begin{alignat*}1
|\phi_i(\vx)-\phi_i(\vy)|\leq L|\vx-\vy| \text{ for } \vx,\vy \in [0,1]^{\Da-1}, i=1,...,M 
\end{alignat*}
such that $S$ is covered by the images
of the maps $\phi_i$.
\end{definition}
We call $L$ a Lipschitz constant for the maps $\phi_i$. By definition the empty set
lies in Lip$(\Da,\M,L)$ for any positive integers $\M$ and $\Da>1$ and any 
non-negative $L$.
\begin{definition}[Adelic-Lipschitz system]\label{defALS}
An adelic-Lipschitz system ($\ALS$) 
$\en_K$ on $K$ (of dimension $n$) is
a set of continuous maps
\begin{alignat*}1
N_v: K_v^{n+1}\rightarrow \Gamma_v \quad v \in M_K
\end{alignat*}
such that for $v \in M_K$ we have
\begin{alignat*}3
(i)&\text{ } N_v({\vz })=0 \text{ if and only if } {\vz} ={\vnull},\\
(ii)&\text{ } N_v(\omega {\vz})=|\omega|_v N_v({\vz}) \text{ for all
$\omega$ in $K_v$ and all ${\vz}$ in $K_v^{n+1}$},\\
(iii)&\text{ if $v \mid \infty$: }\{{\vz};N_v({\vz})=1\} \text{
is in $Lip(d_v(n+1),\M_v,L_v)$ for some $\M_v, L_v$},\\
(iv)&\text{ if $v \nmid \infty$: }N_v({\vz_1}+{\vz_2})
\leq \max\{N_v({\vz}_1),N_v({\vz}_2)\} \text{ for all 
${\vz}_1,{\vz}_2$ 
in $K_v^{n+1}$}.
\end{alignat*}
Moreover we assume that 
\begin{alignat}3
\label{Nvmaxnorm}
N_v(\vz)=\max\{|z_0|_v,...,|z_n|_v\}
\end{alignat}
for all but a finite number
of $v \in M_K$. 
\end{definition}
To deduce our results we will use an $\ALS$ with (\ref{Nvmaxnorm}) for all finite places $v$.
This simplifies the notation and arguments in the sequal considerably. Therefore we assume from now on 
\begin{alignat}3
\label{Nfinmaxnorm}
N_v(\vz)=\max\{|z_0|_v,...,|z_n|_v\} \text{ \qquad for all $v\nmid \infty$.}
\end{alignat}
So the functions $N_v$ with $v\nmid \infty$ are as in Masser and Vaaler's \cite{1}
and the subset of $N_v$ with $v\mid \infty$ 
defines an $(r,s)$-Lipschitz system (of dimension $n$)
in the sense of \cite{1}. However, contrary to Masser and Vaaler we will have to define a uniform $\ALS$
on the collection of all number fields of degree $\e$, as introduced in \cite{art2}. Therefore 
we will use the terminology of \cite{art2}.
With $\M_v$ and $L_v$ from $(iii)$ we define
\begin{alignat*}1
\M_{\en_K}&=\max_{v\mid \infty}\M_v,\\
L_{\en_K}&=\max_{v\mid \infty}L_v.
\end{alignat*}
The set defined in $(iii)$ is the boundary 
of the set ${\bf B}_v=\{{\vz};N_v({\vz})<1\}$
and therefore ${\bf B}_v$ is a bounded symmetric open star-body
in $\IR^{n+1}$ or $\IC^{n+1}$ (see also \cite[p.431]{1}). In particular ${\bf B}_v$ has a finite volume $V_v$.\\

Let us consider the system
where $N_v$ is as in (\ref{Nvmaxnorm}) for all places $v$.
If $v$ is an infinite place then  
${\bf B}_v$ is a
cube for $d_v=1$ and the complex analogue
if $d_v=2$. Their boundaries are clearly
in Lip$(d_v(n+1),M_v,L_v)$ most naturally
with $M_v=2n+2$ maps and $L_v=2$
if $d_v=1$ and 
with $M_v=n+1$ maps and for example $L_v=2\pi\sqrt{2n+1}$
if $d_v=2$.
This system is called the standard
adelic-Lipschitz system.\\

We return to general adelic-Lipschitz systems.
We claim that for any $v\in M_K$ there is a $c_v$ in the value group
$\Gamma_v^*=\Gamma_v\backslash\{0\}$ with
\begin{alignat}3
\label{Nineq1}
N_v({\vz})\geq c_v\max\{|z_0|_v,...,|z_n|_v\}
\end{alignat}
for all $\vz=(z_0,...,z_n)$ in $K_v^{n+1}$.
For if  $v$ is archimedean then ${\bf B}_v$ is
bounded open and it contains the origin.
Since $\Gamma_v^*$ contains arbitrary small
positive numbers the
claim follows by $(ii)$.
Now for $v$ non-archimedean it is trivially true by (\ref{Nfinmaxnorm}) and we can choose $c_v=1$.\\

So let $\en_K$ be an $\ALS$ on $K$ of dimension $n$. For
every $v$ in $M_K$ let
$c_v$ be an element of $\Gamma_v^*$,
such that $c_v\leq 1$ and (\ref{Nineq1}) holds.
Recall we can assume $c_v=1$ for all finite places $v$.
We define
\begin{alignat}3
\label{defcfin}
C^{fin}_{\en_K}&=\prod_{v\nmid \infty}c_v^{-\frac{d_v}{\e}}=1
\end{alignat}
and
\begin{alignat*}3
C^{inf}_{\en_K}&=\max_{v\mid \infty}\{c_v^{-1}\}\geq 1.
\end{alignat*}
Multiplying the finite and the infinite part 
gives rise to another constant 
\begin{alignat}3
\label{defc}
C_{\en_K}&=C^{fin}_{\en_K}C^{inf}_{\en_K}.
\end{alignat}
Besides $\M_{\en_K}$ and $L_{\en_K}$ this is another important quantity for an $\ALS$. We say that \it $\en_K$
is an $\ALS$ with associated constants $C_{\en_K},\M_{\en_K},L_{\en_K}$.\rm \\

In \cite{art1} and \cite{art2} we introduced for an $\ALS$ $\en_K$ on $K$ (of dimension $n$) the quantity $V_{\en_K}^{fin}$. This
quantity depends only on the functions $N_v$ with $v\nmid \infty$ and we have shown in \cite{art1} (first paragraph on p.11) and also in  \cite{art2} (just after equation (3.5)) that if (\ref{Nfinmaxnorm}) holds then
$V_{\en_K}^{fin}=1$. Hence we define
\begin{alignat}3
\label{defVfin1}
V_{\en_K}^{fin}=1.
\end{alignat}
The infinite part is defined by 
\begin{alignat*}3
V_{\en_K}^{inf}=\prod_{v \mid \infty}V_{v}.
\end{alignat*}
By virtue of (\ref{Nineq1}) we observe that
\begin{alignat*}3
V_{\en_K}^{inf}=\prod_{v|\infty} V_v\leq 
\prod_{v|\infty}(2 C^{inf}_{\en_K})^{d_v(n+1)}=
(2 C^{inf}_{\en_K})^{\e(n+1)}.
\end{alignat*} 
We multiply the finite and the infinite part
to get a global volume 
\begin{alignat}3
\label{defVen}
V_{\en_K}=V_{\en_K}^{inf}V_{\en_K}^{fin}.
\end{alignat}
Note that from (\ref{defcfin}), (\ref{defc}), (\ref{defVfin1}) and (\ref{defVen}) we derive
\begin{alignat}3
\label{Venupperbound}
V_{\en_K}\leq (2 C^{inf}_{\en_K}C_{\en_K}^{fin})^{\e(n+1)}=(2 C_{\en_K})^{\e(n+1)}.
\end{alignat}

\subsection{Adelic-Lipschitz heights on $\mathbb{P}^n(K)$}\label{2subsec2}
Let $\en_K$ be an $\ALS$ on $K$
of dimension $n$. Write $\sigma_v$ for the canonical embedding of $K$ into $K_v$, extended componentwise to $K^{n+1}$. Then the height $\henK$ 
on $K^{n+1}$ is defined by
\begin{alignat*}3
\henK(\balf)=\prod_{v \in M_K} N_v(\sigma_v(\balf))^{\frac{d_v}{\e}}.
\end{alignat*}
Thanks to the product formula
and $(ii)$ from Subsection \ref{subsecdefALS}, $\henK(\balf)$ does not change
if we multiply each coordinate of $\balf$ with a fixed element of $K^*$.
Therefore $\henK$ is well-defined on $\IP^n(K)$ by setting
\begin{alignat*}3
\henK(P)=\henK(\balf)
\end{alignat*}
where $P=(\alpha_0:...:\alpha_n) \in \IP^n(K)$ and $\balf=(\alpha_0,...,\alpha_n) \in K^{n+1}$.
Multiplying (\ref{Nineq1}) over all places with 
suitable multiplicities yields 
\begin{alignat}3
\label{HAquiv}
\henK(P)\geq C_{\en_K}^{-1} H(P)
\end{alignat}
for $P\in \IP^n(K)$.

\subsection{Adelic-Lipschitz systems on a collection of number fields}\label{subsec1.1}
We define $\coll_\e$ as the collection of all number fields $K$ of degree $\e$
\begin{alignat*}1
\coll_\e=\{K\subseteq \Qbar;[K:\IQ]=\e\}.
\end{alignat*}
Let $\en$ be a collection of adelic-Lipschitz systems $\en_K$ of dimension $n$ - one for each $K$ of $\coll_\e$. Then we call $\en$ an \em adelic-Lipschitz system $(\ALS)$ on 
$\coll_\e$ of dimension $n$. \rm
We say $\en$ is a \em uniform \rm $\ALS$ on $\coll_\e$
of dimension $n$
with associated constants $C_{\en},\M_{\en},L_{\en}$ in $\IR$ if
the following holds: 
for each $\ALS$ $\en_K$ of the collection $\en$
we can choose associated constants
$C_{\en_K},\M_{\en_K},L_{\en_K}$ satisfying
\begin{alignat*}1
C_{\en_K}\leq C_{\en},\quad
\M_{\en_K}\leq \M_{\en},\quad
L_{\en_K}\leq L_{\en}.
\end{alignat*}
A standard example for a uniform $\ALS$ on $\coll_\e$ (of dimension $n$) is given as 
follows: for each $K$ in $\coll_\e$ choose the standard $\ALS$ on $K$ (of dimension $n$)
so that $N_v$ is as in (\ref{Nvmaxnorm}) for each
$v$ in $M_K$. For this system we may choose
$C_{\en}=1$, $\M_{\en}=2n+2$ and $L_{\en}=2\pi\sqrt{2n+1}$.

\subsection{Adelic-Lipschitz heights on $\mathbb{P}^n(\IQ;\e)$}\label{ALHoncoll}
Let $P=(x_0:...:x_n)\in \IP^n(\Qbar)$ and define $\IQ(P)=\IQ(...,x_i/x_j,...)$ ($0\leq i,j\leq n$; $x_j\neq 0$).
Then we define the degree of $P$ (over $\IQ$) as $[\IQ(P):\IQ]$. Write 
$\IP^n(\IQ;\e)$ for the set of points $P$ in $\IP^n(\Qbar)$ with $[\IQ(P):\IQ]=\e$.
Let $\en$ be an $\ALS$ of dimension $n$ on $\coll_\e$.
Now we can define heights on $\IP^n(\IQ;\e)$.
Let $P\in \IP^n(\IQ;\e)$ so that $\IQ(P)\in \coll_\e$.
According to Subsection \ref{2subsec2} we know that $H_{\en_{K}}(\cdot)$
defines a projective height on $\IP^n(K)$ for each $K$ in $\coll_\e$.
Now we define
\begin{alignat*}1
\hen(P)=H_{\en_{\IQ(P)}}(P).
\end{alignat*}
If $\en$ is the standard adelic-Lipschitz system on $\coll_\e$ as defined in Subsection \ref{subsec1.1}
then $\hen$ is simply the multiplicative Weil height $H$ on $\IP^n(\Qbar)$ (as defined in \cite[p.16]{BG}) restricted to $\IP^n(\IQ;\e)$.

\section{Preliminary results}\label{introchap4}
For $K$ a number field let $\IP^n(K/\IQ)$ be the set of primitive points in $\IP^n(K)$ 
\begin{alignat*}1 
\IP^n(K/\IQ)=\{P\in \IP^n(K);\IQ(P)=K\}.
\end{alignat*}
Let $\en_K$ be an adelic-Lipschitz system of dimension $n$ on $K$.
Then $\henK(\cdot)$ defines a height on $\IP^n(K)$.
Now (\ref{HAquiv}) combined with
Northcott's Theorem implies that the counting function
\begin{alignat*}1 
Z_{\en_K}(\IP^n(K/\IQ),\X)=|\{P\in \IP^n(K/\IQ);\henK(P)\leq \X\}|
\end{alignat*}
is finite for all $\X$ in $[0,\infty)$.
The main result Theorem 3.1 in \cite{art1} gives a precise estimate for this counting function.
Here we need only a special case of Corollary 3.2 in \cite{art1} which by itself is a special case of
Theorem 3.1 in \cite{art1}.
Recall the definitions of $S_K(n)$ from (\ref{Schanuelconst}) and  $V_{\en_K}$ from (\ref{defVen}).
\begin{theorem}\label{prop2}
Let $K$ be a number field of degree $\e$. Let $\en_K$ be an adelic-Lipschitz system of dimension $n$ 
on $K$ with associated constants $C_{\en_K},L_{\en_K},\M_{\en_K}$
and write
\begin{alignat*}3
A_{\en_K}&=\M_{\en_K}^{\e}(C_{\en_K}(L_{\en_K}+1))^{\e(n+1)-1}.
\end{alignat*}
Then as $\X>0$ tends to infinity we have
\begin{alignat*}3
Z_{\en_K}(\IP^n(K/\IQ),\X)=
&2^{-r_K(n+1)}\pi^{-s_K(n+1)}V_{\en_K}S_K(n)\X^{\e(n+1)}\\
+&O(A_{\en_K} R_K h_K\delta(K)^{-\e(n+1)/2+1}\X^{\e(n+1)-1}\mathfrak{L}_0)
\end{alignat*}
where 
\begin{alignat*}3
\mathfrak{L}_0&=\log\max\{2,2C_{\en_K}\X\} \text{ if }(n,\e)=(1,1)\text{ and }\mathfrak{L}_0=1 \text{ otherwise}
\end{alignat*}
and the implied constant in the $O$ depends only on $n$ and $\e$.
\end{theorem}
Now let $\en$ be a uniform $\ALS$ on $\coll_\e$ of dimension $n$.
Then $\hen(\cdot)$ defines a height on $\IP^n(\IQ;\e)$ and (\ref{HAquiv}) implies for any $P\in \IP^n(\IQ;\e)$
\begin{alignat*}3
\hen(P)\geq C_{\en}^{-1} H(P).
\end{alignat*}
Again by Northcott's Theorem we conclude
that the associated counting function $Z_{\en}(\IP^n(\IQ;e),\X)$ (which denotes the number of points 
$P$ in $\IP^n(\IQ;\e)$ with $\hen(P)\leq \X$) is finite for all $\X$ in $[0,\infty)$. 
Bearing in mind the definitions of $S_K(n)$ and  $V_{\en_K}$ from (\ref{Schanuelconst}) and (\ref{defVen}) we define the sum
\begin{alignat}1
\label{kkonst}
\Ce_{\en}=\Ce_{\en}(\IQ,\e,n)=\sum_{K\in \coll_\e}2^{-r_K(n+1)}\pi^{-s_K(n+1)}V_{\en_K}S_K(n).
\end{alignat}
We claim that the sum in (\ref{kkonst}) 
converges if $n$ is large enough.
Now we can state the main result of \cite{art2}. Again we need only a simpler form and so we 
state only this special case of the result.
\begin{theorem}\label{main theorem}
Let $\e,n$ be positive integers. Suppose $\en$ is a uniform adelic-Lipschitz system of dimension $n$
on $\coll_\e$, the collection of all number fields of degree $\e$, with associated constants $C_{\en},\M_{\en}$ and $L_{\en}$. Write
\begin{alignat*}1
A_{\en}&=\M_{\en}^{\e}(C_{\en}(L_{\en}+1))^{\e(n+1)-1}.
\end{alignat*}
Suppose that either $\e=1$ or
\begin{alignat*}1
n>{5\e}/{2}+4+2/\e.
\end{alignat*}
Then the sum in (\ref{kkonst}) converges and 
as $\X>0$ tends to infinity we have
\begin{alignat*}1
Z_{\en}(\IP^n(\IQ;\e),\X)=\Ce_{\en}\X^{\e(n+1)}
+O(A_{\en}\X^{\e(n+1)-1}\mathfrak{L}_0),
\end{alignat*}
where $\mathfrak{L}_0=\log\max\{2,2C_{\en}\X\}$ if
$(\e,n)=(1,1)$ and $\mathfrak{L}_0=1$ otherwise.
The constant in $O$ depends only on
$\e$ and $n$.
\end{theorem}
The following upper bounds are immediate consequences of Schmidt's Theorem in \cite{22}.
\begin{lemma}\label{Zenuppbounds}
Suppose $\en_K$ is an adelic-Lipschitz system (of dimension $n$) on $K$ with associated constants $C_{\en_K}, M_{\en_K}, L_{\en_K}$.
Then 
\begin{alignat}3
\label{Zenupperbound1}
Z_{\en_K}(\IP^n(K),\X)\leq c_1 (C_{\en_K}\X)^{\e(n+1)}.
\end{alignat}
One can choose $c_1=2^{\e(n+4)+n^2+10n+11}$.\\
Now suppose $\en$ is a uniform adelic-Lipschitz system (of dimension $n$) on $\coll_\e$ with associated constants $C_{\en}, M_{\en}, L_{\en}$.
Then 
\begin{alignat}3
\label{Zenupperbound2}
Z_{\en}(\IP^n(\IQ;\e),\X)\leq  c_2(C_{\en}\X)^{\e (\e+n)}.
\end{alignat}
Here one can choose $c_2=2^{\e(\e+n+3)+\e^2+n^2+10\e+10n}$.
\end{lemma}
\begin{rproof}
By (\ref{HAquiv}) we know $\henK(P)\geq C_{\en_K}^{-1}H(P)$ for $P\in \IP^n(K)$, and similar for 
$P\in \IP^n(\IQ,\e)$ one has $\hen(P)\geq C_{\en}^{-1}H(P)$. Thus the statements follow from inequality (1.4) in \cite[Theorem]{22}.
\end{rproof}
We will also use Vinogradov's notation $A\ll B$ (or equivalently $B\gg A$) meaning that there exists a positive constant $c$
depending solely on $\e$ and $n$ (unless specified otherwise) such that $A\leq c B$.
We remind the reader to the definition of the invariant $\delta(K)=\inf\{H(\alpha);K=\IQ(\alpha)\}$.
The following arguments will be used several times. It is therefore convenient to state them as two individual lemmas.
\begin{lemma}\label{SilvermanPrimitiveelement}
Let $K$ be a number field of degree $e>1$ and let $P\in \IP^n(K)$ with $\IQ(P)=K$. Then
\begin{alignat*}3
H(P)&\geq \frac{1}{\e(n+1)}\delta(K),\\
\delta(K)&\geq \e^{-\frac{1}{2(\e-1)}}|\Delta_K|^{\frac{1}{2\e(\e-1)}}.
\end{alignat*}
\end{lemma}
\begin{rproof}
Let us start with the first inequality.
Let $P=(\alpha_0:...:\alpha_n)$ then we can assume that one of the coordinates of $P$
is $1$. Hence $K=\IQ(\alpha_0,...,\alpha_n)$. Now Lemma 3.3 in \cite{art1} gives an element 
$\alpha=\sum_{i=0}^n m_i\alpha_i$ with $0\leq m_i<\e$ in $\IZ$ and $K=\IQ(\alpha)$.
Therefore $H(\alpha)\geq \delta(K)$, and a straightforward computation shows that $H(\alpha)\leq \e(n+1)H(P)$.
This proves the first inequality.
The second inequality is a a special case of Silverman's inequality (\cite[Theorem 2]{9}), but see also
(4.10) and (4.12) in \cite{art2} (with $k=\IQ$ and $m=1$) for more details.
\end{rproof}
\begin{lemma}\label{dyadicsum}
Let $\eta$ be a real number satisfying $\eta<-\e(\e+1)$.
Then we have
\begin{alignat*}3
\sum_{K\in \coll_\e}\delta(K)^{\eta}&\ll_{\eta} 1.
\end{alignat*}
\end{lemma}
\begin{rproof}
This lemma is an immediate consequence of Lemma 4.1 and Lemma 4.3 in \cite{art2}.
\end{rproof}

\section{Reformulation of Theorem \ref{th1} step two: choosing the right Adelic Lipschitz system}

Let $M$ be the Mahler measure on polynomials in one variable with complex coefficients as in \cite{1}.
For each number field $\L$ we define an $ALS$ (of dimension $n$) denoted by $\enML$ by choosing
\begin{alignat}2
\nonumber
&N_v(z_0,...,z_n)=M(z_0x^n+\cdots +z_n)\quad &(v \mid \infty),\\
\label{ALS2}
&N_v(z_0,...,z_n)=\max\{|z_0|_v,...,|z_n|_v\}\quad &(v \nmid \infty).
\end{alignat}
Here $v$ runs over all places in $M_\L$. Masser und Vaaler have shown that $M$ satisfies
$(i),(ii),(iii)$ from Definition \ref{defALS} and with $N_v$ as in (\ref{ALS2}) clearly $(iv)$ is satisfied as well.
Therefore $\henML$ defines an adelic-Lipschitz height height on $\IP^n(\L)$.
Now $M_v$ and $L_v$ depend on $v$ (and $n$), but more precisely they depend only on $d_v\in\{1,2\}$ (and $n$).
Hence $\M_{\enML}$ and $L_{\enML}$ can be chosen independently of $\L$, depending solely on $n$.
Recall the definition of $c_v$ from (\ref{Nineq1}) in Section \ref{subsecdefALS}.
For $v \nmid \infty$ we have $c_v=1$
and for $v | \infty$ we may use $c_v=2^{-n}$
(see \cite[Lemma 2.2, p.56]{3}).
Hence we may set 
\begin{alignat*}1
C_{\enML}=2^n.
\end{alignat*}
So we have shown that we can choose associated constants $C_{\enML}=2^n$, $\M_{\enML}$ and $L_{\enML}$ of the adelic-Lipschitz system $\enML$ 
depending only on $n$.\\

Now let $K$ run over all fields in $\coll_\e$.
The collection of adelic-Lipschitz systems $\enMK$, one
for each number field in $\coll_\e$, defines an adelic-Lipschitz system denoted by $\enM$ on $\coll_\e$.
Then the corresponding height $\henM$ is defined on $\IP^n(\IQ;\e)$.
Furthermore we just have seen that the associated constants $C_{\enMK}=2^n, \M_{\enMK}, L_{\enMK}$ of $\enMK$ may be chosen uniformly,
depending solely on $n$. Thus $\enM$ defines a uniform $\ALS$ on $\coll_\e$ with associated constants $C_{\enM}=2^n, \M_{\enM}, L_{\enM}$.\\

The proofs of our results require also the analogous heights to $\henMK$ and $\henM$ on $\IP^{n}$ but with $n$
replaced by smaller values. By abuse of notation we will use the same symbols
$\henMK$ and $\henM$ for the analogous heights on e.g. $\IP^{n-1}$. But this will cause no confusion.\\

We have a one-to-one correspondence
between monic polynomials in $K[x]$ of degree not exceeding  $n$ and $\IP^n(K)$
\begin{alignat*}1
f_0x^n+\cdots +f_{n-1}x+f_n\longleftrightarrow (f_0:...:f_n).
\end{alignat*}
In this way $\henMK$ can be considered as a function on the monic polynomials in $K[x]$ of degree $\leq n$. In this case we will use 
$M_0$ instead of $\henMK$, so that 
$M_0(f)=\henMK(P_f)$, where $P_f=(f_0:...:f_n)$ and 
$f=f_0x^n+\cdots +f_n$. However, we have also to count monic polynomials whose coefficents do not lie in $K$. Therefore it is convenient to notice that  $M_0$ provides a definition on non-zero polynomials in $\Qbar[x]$ of degree at most $n$. This can be seen in the following way; if $\L$ is any number field containing the coefficients of the non-zero polynomial $f=\alpha_0x^n+\cdots +\alpha_n$ then we set
\begin{alignat*}1
M_0(f)=\henML(P_f)=\prod_{v\in M_\L}N_v(\sigma_v(\alpha_0),...,\sigma_v(\alpha_n))^{d_v/[\L:\IQ]}.
\end{alignat*}
But just as for the usual Weil height it is easy to see
that this definition does not depend on the field $\L$ containing the coordinates and thus $M_0$ is well-defined on the
non-zero polynomials in $\Qbar[x]$ of degree at most $n$.
The Mahler measure $M$ is multiplicative which
together with Gauss' Lemma implies
\begin{alignat}1\label{multiplicativity}
M_0(gh)=M_0(g)M_0(h)
\end{alignat}
for $g,h$ in $\Qbar[x]\backslash 0$ with $\deg gh\leq n$.\\

In the next section we shall see that the proofs of all the theorems can essentially be reduced to finding (asymptotic) estimates for $Z_{\enM}(\IP^n(\IQ;\e),\X)$ as given in Theorem \ref{main theorem}.

\section{Proofs of the Theorems}
We remind the reader that $K$ denotes a number field of degree $\e$.
As mentioned in the introduction for $\e=1$ or $n=1$ all our theorems are covered by results of
Schmidt \cite{22}, Masser and Vaaler \cite{37}, \cite{1} and the author \cite{art1}. From now on we assume 
\begin{alignat*}1
\e>1 \text{ and } n>1.
\end{alignat*}
We start with the set
\begin{alignat*}1 
\grmn=\{f\in K[x];f\text{ monic,} \deg f\leq n, \IQ(P_f)=K,
M_0(f)\leq \X\}.
\end{alignat*}
Recall that $\IP^n(K/\IQ)$ is the set of primitive points in $\IP^n(K)$ and
$Z_{\enMK}(\IP^n(K/\IQ),\X)$ is its counting function with respect to $\henMK$.
Then of course
\begin{alignat}1 
\label{grmcountfct}
|\grmn|=Z_{\enMK}(\IP^n(K/\IQ),\X).
\end{alignat}
For any $f$ in $\grmn$ one has
\begin{alignat*}1 
\X\geq M_0(f)=\henMK(P_f).
\end{alignat*}
Moreover we know  
$\henMK(P_f)\geq C_{\enMK}^{-1}H(P_f)=2^{-n}H(P_f)$.
Now $f\in \grmn$ implies $K=\IQ(P_f)$ and hence we can apply Lemma \ref{SilvermanPrimitiveelement} to deduce 
\begin{alignat*}1 
H(P_f)\geq \frac{1}{\e(n+1)}\delta(K).
\end{alignat*}
Note also that the Mahler measure of a monic polynomial is at least $1$ and therefore $M_0(f)\geq 1$. So whenever $\grmn$ is non-empty we have
\begin{alignat}1 
\label{heightcomp1'}
&\X\geq 1,\\
\label{heightcomp1}
&\X\geq \frac{\delta(K)}{2^{n}\e(n+1)}\gg \delta(K).
\end{alignat}
For a subfield $k$ of $K$ let $\Hom_{k}(K)$ be the set of $k$-invariant 
field homomorphisms from $K$ to its Galois closure $K_G$ over $\IQ$.
\newline
Let $\Mcp(n,\X)$ be the set of all monic, irreducible 
polynomials $f$ of degree $n$ in $K[x]$, with  
$\sigma f$ are pairwise coprime as $\sigma$ runs over $\Hom_{\IQ}(K)$ 
and $M_0(f)\leq \X$. Here the homomorphisms $\sigma$ act on the coefficients of the polynomials.
Note that the coprimality of the polynomials $\sigma f$
implies $\IQ(P_f)=K$. Hence
\begin{alignat*}1
\Mcp(n,\X)=\{& f\in \grmn\backslash \grmnm1;
f\text{ irreducible over $K$,}\\
\nonumber&\text{$\sigma f$ pairwise coprime ($\sigma \in \Hom_{\IQ}(K)$)}\}.
\end{alignat*}
\begin{lemma} 
\label{lemma1chap5}
We have
\begin{alignat}1\label{eqlemma1chap5}
 Z_K(\e,n,X)=n |\Mcp(n,X^n)|.
\end{alignat}
\end{lemma}
\begin{rproof}
We will show that the map that sends $\beta$ to its monic minimal polynomial over $K$
defines a $n$-to-one correspondence between the set $S_K(\e,n,X)=\{\beta\in \Qbar;[\IQ(\beta):\IQ]=\e n, [K(\beta):K]=n,
H(\beta)\leq X\}$ (corresponding to the counting function $Z_K(\e,n,X)$) and the set $\Mcp(n,X^n)$.\\
Let $f$ be in $K[x]$ irreducible with $\deg f=n$.
Then $f$ has $n$ zeros, they are pairwise
distinct and, of course, each of them has degree $n$ over $K$. Therefore we get a factor $n$. On the other hand every $\beta$ with $[K(\beta):K]=n$
is a zero of exactly one irreducible monic polynomial
$f$ in $K[x]$. We factor  
$f=(x-\beta_1)\cdots (x-\beta_n)$. Then
\begin{alignat*}1
M_0(f)=M_0(x-\beta_1)\cdots M_0(x-\beta_n).
\end{alignat*}
Since $f$ is irreducible all the zeros of
$f$ have the same height. But $H(\alpha)=M_0(x-\alpha)$ for any $\alpha \in \Qbar$ and
so we get
\begin{alignat}1
\label{M0H}
M_0(f)=H(\beta_1)^n.
\end{alignat}
This explains the power $X^n$.\\
Now let $D_{\beta,\IQ}$ be the monic minimal polynomial of $\beta$ over $\IQ$.
Then clearly $f|D_{\beta,\IQ}$. If the $\sigma f$ are not pairwise coprime then
\begin{alignat*}1
\prod_{\Hom_{\IQ}(K)} \sigma f,
\end{alignat*}
which of course lies in $\IQ[x]\backslash \IQ$, cannot be irreducible over $\IQ$.
Hence $[\IQ(\beta):\IQ]<|\Hom_{\IQ}(K)|\deg f=\e n$ which means $\beta \notin S_K(\e,n,X)$.
Next we notice that 
for any $\sigma$ of $\Hom_{\IQ}(K)$ we have
\begin{alignat*}1
\sigma f|\sigma D_{\beta,\IQ}=D_{\beta,\IQ}|\prod_{\Hom_{\IQ}(K)} \sigma f.
\end{alignat*}
Now suppose the $\sigma f$ are pairwise coprime then
\begin{alignat*}1
\prod_{\Hom_{\IQ}(K)} \sigma f|D_{\beta,\IQ}
\end{alignat*}
and we end up with $[\IQ(\beta):\IQ]=|\Hom_{\IQ}(K)|\deg f=\e n$ which shows $\beta \in S_K(\e,n,X)$.
This completes the proof.
\end{rproof}
To count $|\Mcp(n,\X)|$ via $|\grmn|$ another
two sets are required.
First we define the subset
\begin{alignat*}1
\Mred(n,\X)=\{f\in \grmn\backslash \grmnm1;f\text{ reducible over $K$}\}.
\end{alignat*}
So $\Mred(n,\X)$ is the set of all monic 
reducible polynomials $f$ of degree $n$ in $K[x]$ with
$K=\IQ(P_f)$ and $M_0(f)\leq \X$. Finally let
\begin{alignat*}1
\Mncp(n,\X)=\{& f \in \grmn\backslash \grmnm1;
f\text{ irreducible over $K$,}\\
&\sigma f\text{ not pairwise coprime } 
(\sigma \in \Hom_{\IQ}(K))\}.
\end{alignat*}
Immediately from the definition we get
\begin{alignat}1
\label{Zer}
\Mcp(n,\X)=\grmn\backslash\left(\grmnm1\cup \Mred(n,\X)\cup \Mncp(n,\X)\right).
\end{alignat}
In particular
\begin{alignat}1
\label{Zerest}
|\Mcp(n,\X)|\leq |\grmn|.
\end{alignat}
From (\ref{grmcountfct}) we get
\begin{alignat}1
\label{Zer2}
\sum_{K\in \coll_\e} |\grmn|=
\sum_{K\in \coll_\e}Z_{\enMK}(\IP^n(K/\IQ),\X)=
Z_{\enM}(\IP^n(\IQ;\e),\X).
\end{alignat}
Now (\ref{ZmnXub2}) and Lemma \ref{lemma1chap5}
yields
\begin{alignat*}2
Z(\e,n,X)\leq \sum_{K\in \coll_\e}Z_K(\e,n,X)=n\sum_{K\in \coll_\e}|\Mcp(n,X^n)|.
\end{alignat*}
Taking into account (\ref{Zerest}) and (\ref{Zer2}) gives
\begin{alignat}2
\label{cha4red2'}
Z(\e,n,X)\leq
n Z_{\enM}(\IP^n(\IQ;\e),X^n).
\end{alignat}
In order to obtain asymptotic estimates more care is needed. 
Combining (\ref{3gl3}), (\ref{eqlemma1chap5}) and (\ref{Zer}) we get as $X>0$ tends to infinity
\begin{alignat*}2
Z(\e,n,X)=n \sum_{K\in \coll_\e}|\grmnnew|
+&O(\sum_{K\in \coll_\e}|\grmnmnew1|)\\
+&O(\sum_{K\in \coll_\e}|\Mred(n,X^n)|)\\
+&O(\sum_{K\in \coll_\e}|\Mncp(n,X^n)|)\\
+&O(\sum_{\f|n \atop{1<\f\leq \e}}\sum_{\L \in \coll_{\f\e}}|\MFcp(n/\f,X^{n/\f})|).
\end{alignat*}
Applying (\ref{Zerest}) gives $|\MFcp(n/\f,X^{n/\f})|\leq |\MFall(n/\f,X^{n/\f})|$ and then applying (\ref{Zer2}) for the first, second and the last term yields
\begin{alignat}2
\label{cha4red2}
Z(\e,n,X)=n Z_{\enM}(\IP^n(\IQ;\e),X^n)
+&O(Z_{\enM}(\IP^{n-1}(\IQ;\e),X^n))\\
\label{2cha4red2}
+&O(\sum_{K\in \coll_\e}|\Mred(n,X^n)|)\\
\label{3cha4red2}
+&O(\sum_{K\in \coll_\e}|\Mncp(n,X^n)|)\\
\label{4cha4red2}
+&O(\sum_{\f|n \atop{1<\f\leq \e}}Z_{\enM}(\IP^{n/\f}(\IQ;\f\e),X^{n/\f})).
\end{alignat}
To handle the error terms we need good uniform upper bounds for $Z_{\enML}(\IP^n(\L),\X)$ and $Z_{\enML}(\IP^n(\L/\IQ),\X)$.
\begin{lemma} 
\label{lemmauppbounds2}
Let $\L$ be a number field and let $m\leq n$ be a positive integer. Then
\begin{alignat}1
\label{S1'}
Z_{\enML}(\IP^m(\L),\X)\ll_{[\L:\IQ]} \X^{[\L:\IQ](m+1)}.
\end{alignat}
\end{lemma}
\begin{rproof}
Recall that $C_{\enML}=2^m$ and $m\leq n$. Thus the statement follows from (\ref{Zenupperbound1}) in Lemma \ref{Zenuppbounds}.
\end{rproof}
\begin{lemma}\label{lemmauppbounds}
Let $\L$ be a number field and let $m\leq n$ be a positive integer. Then
\begin{alignat}1
\label{S1}
Z_{\enML}(\IP^m(\L/\IQ),\X)\ll_{[\L:\IQ]} 
\frac{R_\L h_\L}{\delta(\L)^{\frac{[\L:\IQ](m+1)}{2}}}\X^{[\L:\IQ](m+1)}.
\end{alignat}
\end{lemma}
\begin{rproof}
The case $\L=\IQ$ is covered by the preceeding lemma, so we can assume $[\L:\IQ]>1$.
If $Z_{\enML}(\IP^m(\L/\IQ),\X)=0$ then the claim is certainly true. Now assume $Z_{\enML}(\IP^m(\L/\IQ),\X)>0$.
In this case we know from (\ref{grmcountfct}) and (\ref{heightcomp1}) that  $\X\gg_{[\L:\IQ],m} \delta(\L)$.
For $[\L:\IQ]>1$ Theorem \ref{prop2} immediately implies 
\begin{alignat*}1
Z_{\enML}(\IP^m(\L/\IQ),\X)\ll_{[\L:\IQ],m,C_{\enML},M_{\enML},L_{\enML}} 
&\frac{R_\L h_\L}{|\Delta_\L|^{\frac{(m+1)}{2}}}V_{\enML}\X^{[\L:\IQ](m+1)}\\
&+\frac{R_\L h_\L}{\delta(\L)^{\frac{[\L:\IQ](m+1)}{2}-1}}\X^{[\L:\IQ](m+1)-1}.
\end{alignat*}
Recall that $C_{\enML},M_{\enML},L_{\enML}$ depend only on $m$; but $m\leq n$ and thus they are $\ll 1$.
Therefore and due to (\ref{Venupperbound}) we have $V_{\enML}\ll 1$. Moreover we get $\X\gg_{[\L:\IQ]} \delta(\L)$ and hence 
\begin{alignat*}1
Z_{\enML}(\IP^m(\L/\IQ),\X)\ll_{[\L:\IQ]} 
\frac{R_\L h_\L}{|\Delta_\L|^{\frac{(m+1)}{2}}}\X^{[\L:\IQ](m+1)}+
\frac{R_\L h_\L}{\delta(\L)^{\frac{[\L:\IQ](m+1)}{2}}}\X^{[\L:\IQ](m+1)}.
\end{alignat*}
Now Lemma 4.5 in \cite{art2} gives $|\Delta_\L|\gg_{[\L:\IQ]} \delta(\L)^{[\L:\IQ]}$. This proves the lemma.
\end{rproof}
Note that by Siegel-Brauer's Theorem $R_Kh_K\ll|\Delta_K|^{1/2+1/(40\e(\e-1))}$ and recall the inequality 
$\delta(K)\gg |\Delta_K|^{\frac{1}{2\e(\e-1)}}$ from Lemma \ref{SilvermanPrimitiveelement}. Thus we get 
\begin{alignat}1
\label{deltaRh}
R_Kh_K\ll \delta(K)^{\e(\e-1)+1/20}.
\end{alignat}

\subsection{An upper bound for $|\Mred(n,\X)|$}
In this subsection we will prove an upper bound for the number of polynomials $f\in \Mall(n,\X)$ of degree $n$
that are reducible over $K$. Recall that by definition $\delta(K)\geq 1$ and by (\ref{heightcomp1'}) and (\ref{heightcomp1}) we can assume
$\X\geq 1$ and $\X/\delta(K)\gg 1$.\\

Suppose $f$ factors as
\begin{alignat*}1
f=gh
\end{alignat*}
where $g,h$ are in $K[x]\backslash K$ and monic.
Since $K=\IQ(P_f)\subseteq \IQ(P_g,P_h)\subseteq K$
three cases may occur.
\begin{alignat*}1
(A):&\quad \IQ(P_g)=K,\quad \IQ(P_h)=K,\\
(B):&\quad \IQ(P_g)\subsetneq K,\quad \IQ(P_h)=K,\\
(C):&\quad \IQ(P_g)\subsetneq K,\quad \IQ(P_h)\subsetneq K.
\end{alignat*}
Let $\deg g=p$ so that $1\leq p \leq n-1$ and $\deg h=n-p$. Assume $M_0(f)\leq \X$.
Now $M_0(f)\geq 1$
and hence there exists a positive integer $i$ such that $2^{i-1}\leq M_0(g)< 2^i$ and then the multiplicativity (\ref{multiplicativity}) of $M_0$
gives $M_0(h)\leq 2^{1-i}\X$. For fixed $i$ we will estimate the number of polynomials $f=gh$ in each of the three cases $(A)$, $(B)$
and $(C)$ separately and then we sum over all possible values for $i$, i.e. $i=1,...,[\log_2 \X]+1$.
To simplify the notation we abbreviate $\delta(K)$ to $\delta$.\\

We start with the case $(A)$.
Here we can assume by symmetry that $p\leq n/2$.
To bound the number of polynomials $f=gh$ we apply Lemma \ref{lemmauppbounds} with $\L=K$. Thus for fixed $i$ we get the upper bound
\begin{alignat*}1
&\ll \left(R_Kh_K\delta^{-\frac{\e}{2}(p+1)}(2^i)^{\e(p+1)}\right)\left(R_Kh_K\delta^{-\frac{\e}{2}(n-p+1)}(2^{1-i}\X)^{\e(n-p+1)}\right)\\
&= 2^{\e(n-p+1)}(2^i)^{\e(2p-n)}(R_Kh_K)^2\delta^{-\frac{\e}{2}(n+2)}\X^{\e(n-p+1)}
\end{alignat*}
for the number of $f$. Now if $p<n/2$ then $\sum_i (2^i)^{\e(2p-n)}\ll 1$ where the sum runs over all values $i=1,...,[\log_2 \X]+1$. So in this case we get the upper bound 
\begin{alignat*}1
\ll (R_Kh_K)^2\delta^{-\frac{\e}{2}(n+2)}\X^{\e n}
\end{alignat*}
for the number of polynomials $f=gh$. Now suppose $n=p/2$. Then the sum over $i$ introduces an additional logarithm and we find the upper bound
\begin{alignat*}1
\ll (R_Kh_K)^2\delta^{-\frac{\e}{2}(n+2)}\X^{\e(n/2+1)}\log(\X+2)\ll (R_Kh_K)^2\delta^{-\frac{\e}{2}(n+2)}\X^{\e n}.
\end{alignat*}
Next we use (\ref{deltaRh})
to eliminate $R_Kh_K$. This yields for the number in $(A)$
\begin{alignat*}1
\ll \delta^{-\frac{\e}{2}(n-4\e+6)+0.1}\X^{\e n}.
\end{alignat*}

Next we estimate the number of polynomials in $(B)$. We proceed similar as in $(A)$. But here the situation is not symmetric hence we cannot assume
$p\leq n/2$ and moreover we use (\ref{S1'}) with $\L\subsetneq K$ to bound the number of polynomials $g$.
Note also that there are only $\leq 2^{\e}\ll 1$ possibilities for $\L$. 
For fixed $i$ this yields the upper bound
\begin{alignat*}1
\ll R_Kh_K\delta^{-\frac{\e}{2}(n-p+1)}\X^{\e(n-p+1)}2^{-\frac{i\e}{2}(2n-3p+1)}.
\end{alignat*}
Then summing over $i=1,...,[\log_2 \X]+1$ we obtain $3$ different upper bounds depending on whether $2n-3p+1>0$, $2n-3p+1=0$ or $2n-3p+1<0$.
Finally we use $\X/\delta\gg 1$ and (\ref{deltaRh}) to deduce that also all of these $3$ upper bounds are covered by
\begin{alignat*}1
\ll \delta^{-\frac{\e}{2}(n-4\e+6)+0.1}\X^{\e n}.
\end{alignat*}
We are left with the case $(C)$. Here we use (\ref{S1'}) with $\L\subsetneq K$ to bound the number of polynomials $g$ and $h$. 
And again we use the fact that there are only $\leq 2^{\e}\ll 1$ possibilities for $\L$.
Furthermore, by symmetry, we can assume $p\leq n/2$.
Similar as in $(A)$ we obtain the upper bound 
\begin{alignat*}1
\ll \X^{\frac{\e n}{2}}\ll \X^{\frac{\e n}{2}}(\X/\delta)^{\frac{\e n}{2}}\ll \delta^{-\frac{\e }{2}(n-4\e+6)+0.1}\X^{\e n}.
\end{alignat*}
Again we can multiply the error terms arising
from $(A)$, $(B)$ and $(C)$ with $(\X/\delta)^a$
as long as $a\geq 0$. We choose $a$ such that the exponent
on $\X$ is $\e(n+1)-1$. Hence all three error terms are covered by
\begin{alignat*}1
\ll\delta^{-\frac{\e}{2}(n-4\e+8)+1.1}
\X^{\e(n+1)-1}.
\end{alignat*}
Thus we have proven
\begin{alignat}1
\label{FS1}
|\Mred(n,\X)|\ll \delta(K)^{-\frac{\e}{2}(n-4\e+8)+1.1}
\X^{\e(n+1)-1}. 
\end{alignat}

\subsection{An upper bound for $|\Mncp(n,\X)|$}
As in the previous subsection we can assume
$\X\geq 1$ and $\X/\delta(K)\gg 1$.
Recall that $K_G$ is the Galois closure of $K$ over $\IQ$.
Suppose $f$ is in $\grmn$ and 
irreducible over $K_G$. Hence for all
$\sigma \in \Hom_{\IQ}(K)$ the $\sigma f$ are irreducible in $K_G[x]$ and since $\IQ(P_f)=K$ they are pairwise distinct.
Thus they are pairwise coprime.
It follows
\begin{alignat}1
\label{ink1}
&\Mncp(n,\X)\subseteq \\
\nonumber&\{f \in \grmn\backslash \grmnm1;\text{$f$ irreducible over $K$, $f$ reducible over $K_G$}\}.
\end{alignat}
So let $f$ be as above; that is $f\in K[x]$ monic, irreducible over $K$ but reducible over $K_G$, $\deg f=n$ and $\IQ(P_f)=K$. Let 
\begin{alignat*}1
f=g_1\cdots g_s
\end{alignat*}
be its decomposition into prime factors in $K_G[x]$ 
($g_1,...,g_s$ pairwise distinct, monic) 
and let
\begin{alignat*}1
\L=K(P_{g_1})
\end{alignat*}
be the field, gotten by adjoining the coefficients of $g_1$ to $K$.
\begin{lemma} 
We have
\begin{alignat*}1
f=\prod_{\tau \in \Hom_{K}(\L)}\tau g_1.
\end{alignat*}
\end{lemma}
\begin{rproof}
First notice that
\begin{alignat*}1
\prod_{\tau \in \Hom_{K}(\L)}\tau g_1 \in K[x]. 
\end{alignat*}
For $\tau$ as in the product above we have that $\tau g_1$ divides $\tau f=f$.
Since $\IQ(P_{g_1})=\L$ the $\tau g_1$ are pairwise distinct.
For any such $\tau$ there is a $\sigma$ in $\rm{Gal}(K_G/\IQ)$
with $\tau g_1=\sigma g_1$. But $g_1$ is irreducible in $K_G[x]$ and so the $\sigma g_1$ are all 
irreducible in $K_G[x]$. Thus the $\tau g_1$  
are irreducible pairwise distinct divisors of $f$ in $K_G[x]$ and therefore they are also pairwise
coprime.
This yields $\prod_{\tau \in \Hom_{K}(\L)} \tau g_1$ divides $f$.
On the other hand $\prod_{\tau \in \Hom_{K}(\L)} \tau g_1$ is in $K[x]\backslash K$ and monic.
Since  
$f$ is monic and irreducible over $K$ they are equal.
\end{rproof}
Let $f=(x-\beta_1)\cdots (x-\beta_n)$ be the factorisation in $\Qbar[x]$.
The function $M_0$ is defined on polynomials in $\Qbar[x]$ of degree not larger than $n$ and is  multiplicative.
Therefore $M_0(f)=M_0(x-\beta_1)\cdots M_0(x-\beta_n)$.
Now $f$ is irreducible in $K[x]$ so all the zeros
have the same height or equivalently
$M_0(x-\beta_1)=\cdots=M_0(x-\beta_n)$. In particular
$M_0(g_1)=M_0(\tau g_1)$ 
for all $\tau \in \Hom_{K}(\L)$.
We conclude
\begin{alignat*}1
\X\geq M_0(f)=M_0(g_1)^{[\L:K]}.
\end{alignat*}
To bound the cardinality of the set in (\ref{ink1}) above, we proceed as follows: for any intermediate field $\L$ with 
$K\subsetneq \L\subseteq K_G$ we estimate the number of  monic $g\in \L[x]$ with 
\begin{alignat}1\label{degn}
&\deg g [\L:K]=\deg f=n\\
&M_0(g)\leq \X^{\frac{1}{[\L:K]}}.
\end{alignat}
Then we sum these estimates over all fields $\L$.
Hence we have
\begin{alignat*}1
|\Mncp(n,\X)|\leq \sum_{\L\atop K\subsetneq \L\subseteq K_G}|\{g\in \L[x]; g \text{ monic }, \deg g=\frac{n}{[\L:K]}, M_0(g)\leq \X^{\frac{1}{[\L:K]}}\}|.
\end{alignat*}
Note that of course only fields $\L$ with $[\L:K]\mid n$ contribute to the sum above. Hence we can assume
$[\L:K]\mid n$.
Now clearly 
\begin{alignat*}1
Z_{\enML}(\IP^{\frac{n}{[\L:K]}}(\L),\X^{\frac{1}{[\L:K]}})\geq |\{g\in \L[x]; g \text{ monic, } \deg g=\frac{n}{[\L:K]}, M_0(g)\leq \X^{\frac{1}{[\L:K]}}\}|
\end{alignat*}
and thus
\begin{alignat*}1
|\Mncp(n,\X)|\leq \sum_{\L\atop K\subsetneq \L\subseteq K_G}Z_{\enML}(\IP^{\frac{n}{[\L:K]}}(\L),\X^{\frac{1}{[\L:K]}})
\end{alignat*}
Applying Lemma \ref{lemmauppbounds2}, and not forgetting that by (\ref{degn}) $[\L:\IQ]\ll 1$, yields  
\begin{alignat*}1
Z_{\enML}(\IP^{\frac{n}{[\L:K]}}(\L),\X^{\frac{1}{[\L:K]}})\ll \X^{\frac{[\L:\IQ]}{[\L:K]}\left(\frac{n}{[\L:K]}+1\right)}
= \X^{\e\left(\frac{n}{[\L:K]}+1\right)}\leq \X^{\frac{\e n}{2}+\e}.
\end{alignat*}
The degree of $K_G$ is bounded from above by $e!$.
Therefore the number of intermediate fields $\L$ is bounded from above by $2^{\e!}\ll 1$
and so we end up with 
\begin{alignat*}1
|\Mncp(n,\X)| \ll \X^{\frac{\e n}{2}+\e}.
\end{alignat*}
As in the previous subsection we use (\ref{heightcomp1}) to deduce
\begin{alignat}1
\nonumber |\Mncp(n,\X)|&\ll \delta(K)^{-\frac{\e n}{2}+1}\X^{\e(n+1)-1}\\
\label{FS2}
&\leq \delta(K)^{-\frac{\e}{2}(n-4\e+8)+1.1}\X^{\e(n+1)-1}.
\end{alignat}

\subsection{Proof of Theorem \ref{th4}}
Recall that $\enM$ defines a uniform $\ALS$ with $C_{\enM}=2^n$. So (\ref{Zenupperbound2}) in Lemma \ref{Zenuppbounds}
yields
\begin{alignat*}1
Z_{\enM}(\IP^n(\IQ;\e),\X)\leq c_2(2^n\X)^{\e(n+\e)}
\end{alignat*}
where $c_2$ is defined in Lemma \ref{Zenuppbounds}.
This together with (\ref{cha4red2'}) yields immediately the
following bound
\begin{alignat*}1
Z(\e,n,X)\leq n c_2(2X)^{\e n(n+\e)}
\end{alignat*}
and thereby proves Theorem \ref{th4}.

\subsection{Proof of Theorem \ref{th1}}
Recall the fundamental equality (\ref{cha4red2}). We start with the first term on the right-hand side of (\ref{cha4red2}).
Note that $n>\max\{\e^2+\e,10\}\geq 5\e/2+4+2/\e$ unless $\e=3$. But then $5\e/2+4+2/\e=12+1/6$ and $\e^2+\e=12$ and so
$n>\max\{\e^2+\e,10\}$ implies $n>5\e/2+4+2/\e$ always.  Hence we can apply
Theorem \ref{main theorem} to conclude  
\begin{alignat}1
\label{est3cha4}
n Z_{\enM}(\IP^n(\IQ;\e),X^n)=n\Ce_{\enM}(\IQ,\e,n)X^{\e n(n+1)}
+O(X^{\e n(n+1)-n})
\end{alignat}
where
\begin{alignat}1
\label{eq3cha4}
\Ce_{\enM}(\IQ,\e,n)=\sum_{K\in \coll\e} 2^{-r_K(n+1)}\pi^{-s_K(n+1)}V_{\enMK}S_K(n).
\end{alignat}
From (\ref{defVen}) we recall that $V_{\enMK}=V^{inf}_{\enMK}V^{fin}_{\enMK}$.
The volume $V^{inf}_{\enMK}$ has been computed by Masser und Vaaler in \cite[p.435]{1} (in their notation
$V_{\en}$)
\begin{alignat*}1
V^{inf}_{\enMK}=2^{r_K(n+1)}\pi^{s_K(n+1)}
V_{\IR}(n)^{r_K}V_{\IC}(n)^{s_K}.
\end{alignat*}
By definition (\ref{defVfin1}) we have $V^{fin}_{\enMK}=1$ and hence 
\begin{alignat}1\label{VenMK}
V_{\enMK}=2^{r_K(n+1)}\pi^{s_K(n+1)}
V_{\IR}(n)^{r_K}V_{\IC}(n)^{s_K},
\end{alignat}
supporting our main term.\\
Next we consider the second term on the right-hand side of (\ref{cha4red2}). We could use Theorem \ref{main theorem} again, to get an upper bound 
for $Z_{\enM}(\IP^{n-1}(\IQ;\e),X^n)$. However, it is slightly better to proceed as follows. Clearly
\begin{alignat*}1
Z_{\enM}(\IP^{n-1}(\IQ;\e),X^n)= \sum_{K\in \coll_\e}Z_{\enMK}(\IP^{n-1}(K/\IQ),X^n).
\end{alignat*}
Now from (\ref{S1}) and (\ref{deltaRh}) we find 
\begin{alignat}1
\nonumber Z_{\enMK}(\IP^{n-1}(K/\IQ),X^n)&\ll R_Kh_K\delta(K)^{-\e n/2}X^{\e n(n-1)}\\
\label{ineq620} &\ll\delta(K)^{-\e n/2+\e(\e-1)+0.05}X^{\e n(n-1)}.
\end{alignat}
Next note that $n>\{\e^2+\e,10\}\geq 4\e$. But $n>4\e$ implies $-\e n/2+\e(\e-1)+0.05<-\e(\e+1)$ and so we conclude by virtue of Lemma \ref{dyadicsum}  
\begin{alignat}1
\label{est4cha4}
Z_{\enM}(\IP^{n-1}(\IQ;\e),X^n)\ll X^{\e n(n-1)}\ll X^{\e n(n+1)-n},
\end{alignat}
where in the last inequality we may assume $X\gg 1$ because $\henM(P)\gg 1$ for any $P$ in $\IP^{n-1}(\IQ;\e)$.\\
Now appealing to (\ref{FS1}) and (\ref{FS2}) shows that the remaining terms coming from
(\ref{2cha4red2}) and (\ref{3cha4red2}) are bounded by
\begin{alignat*}1
\ll X^{\e n(n+1)-n} \sum_K\delta(K)^{-\frac{\e}{2}(n-4\e+8)+1.1}.
\end{alignat*}
The latter sum is convergent by virtue of Lemma \ref{dyadicsum} provided
$-\frac{\e}{2}(n-4\e+8)+1.1<-\e(\e+1)$ or equivalently $n>6\e-6+2.2/\e$. But $n>\{\e^2+\e,10\}$ implies $n>6\e-6+2.2/\e$
and so we have proved 
\begin{alignat*}1
\sum_K|\Mred(n,X^n)|+\sum_K|\Mncp(n,X^n)|\ll X^{\e n(n+1)-n}.
\end{alignat*}
To bound the last term in (\ref{4cha4red2}) we apply (\ref{Zenupperbound2}). Recalling $C_{\enM}\ll 1$ we find
\begin{alignat*}1
\sum_{\f|n \atop{1<\f\leq \e}}Z_{\enM}(\IP^{n/\f}(\IQ;\f \e),X^{n/\f})
\ll \sum_{\f|n \atop{1<\f\leq \e}} X^{\e n(\f\e+n/\f)}.
\end{alignat*}
Again we may assume $X\gg 1$ because $\henM(P)\gg 1$.
Now for $2\leq \f\leq \e$ we have $\e n(\f\e+n/\f)\leq \e n(n+1)-n$ provided $n\geq \e^2+\e+1/(\e-1)$.
But by hypothesis we have $n>\{\e^2+\e,10\}$ which implies $n\geq \e^2+\e+1/(\e-1)$.
Hence 
\begin{alignat*}1
\sum_{\f|n \atop{1<\f\leq \e}}Z_{\en}(\IP^{n/\f}(\IQ;\f \e),X^{n/\f})\ll X^{\e n(n+1)-n}. 
\end{alignat*}
This completes the proof of Theorem \ref{th1}.

\subsection{Proof of Theorem \ref{th2}}
Again we start with the equality (\ref{cha4red2}).
Note that the extra condition on $\e$ and $n$ in Theorem \ref{th2}
implies that the sum in (\ref{4cha4red2}) is empty.
In the proof of Theorem \ref{th1} we have seen that the $O$-terms in (\ref{cha4red2}), (\ref{2cha4red2}) and (\ref{3cha4red2})
are bounded from above by $\ll X^{\e n(n+1)-n}$, subject to $n>\max\{5\e/2+4+2/\e,4\e,6\e-6+2.2/\e\}$. 
But $\max\{6\e-6+2.2/\e,10\}\geq \max\{5\e/2+4+2/\e,4\e\}$ and clearly $n>\max\{6\e-6+2.2/\e,10\}$ if and only if $n>\max\{6\e-6,10\}$. Therefore the statement of the theorem follows from
(\ref{est3cha4}) and (\ref{eq3cha4}).

\subsection{Proof of Theorem \ref{th5}}
We claim that
\begin{alignat}1
\label{fundeq}
Z(\e,m,n,X)=\sum_K Z_K(\e m,n,X)
\end{alignat}
where the sum runs over fields $K$ of degree  $\e m$ that contain a subfield of degree $\e$.
Recall that $S_K(\e m,n,X)$ denotes the set counted by $Z_K(\e m,n,X)$
and let $S(\e,m,n,X)$ denote the set counted by $Z(\e,m,n,X)$.\\
First we show ``$\leq$''. Suppose $\beta$ lies in $S(\e,m,n,X)$. Hence there exists a field
$k\subseteq \IQ(\beta)$ and a field $K\subseteq \IQ(\beta)$ with $[k:\IQ]=\e$ and $[K:\IQ]=\e m$.
Suppose $k$ is not contained in $K$. Then $\IQ(\beta)$, which has degree $\e m n$, contains the field compositum of $k$ and $K$
which has degree $\f \e m$ for an $\f$ satisfying $1<\f\leq \e\leq \e m$ and $\f|n$. But the latter contradicts the hypothesis of Theorem \ref{th5}.
Hence each $\beta$ in $S(\e,m,n,X)$ lies in at least one $S_K(\e m,n,X)$.
Now we prove the other inequality ``$\geq$''. Of course each $\beta \in S_K(\e m,n,X)$ lies in $S(\e,m,n,X)$.
Now if $\beta$ lies in $S_K(\e m,n,X)$ and in $S_{K'}(\e m,n,X)$ then $\IQ(\beta)$, which has degree $\e m n$, contains the field compositum of the two different fields $K$ and $K'$ which has 
degree $\f \e m$ for an $\f$ satisfying $1<\f\leq \e m$ and $\f|n$; again this contradicts the hypothesis of Theorem \ref{th5}. 
This proves (\ref{fundeq}).\\
Recalling (\ref{grmcountfct}) and then applying Theorem \ref{prop2} with (\ref{VenMK}) gives: as $X>0$ tends to infinity
\begin{alignat}1\label{th5eq1}
|\grmnnew|&=V_{\IR}(n)^{r_K}V_{\IC}(n)^{s_K}S_K(n)X^{\e m n (n+1)}\\
&+O(R_Kh_K\delta(K)^{-\e m(n+1)/2+1}X^{\e m n (n+1)-n}).
\end{alignat}
And thanks to (\ref{deltaRh}) the error term above is covered by 
\begin{alignat}1\label{th5eq2}
\ll \delta(K)^{-\frac{\e m}{2}(n-2\e m+3)+1.05}X^{\e m n(n+1)-n}.
\end{alignat}
Applying Lemma \ref{dyadicsum} shows that the above error term converge when summed over
$\coll_{\e m}$ and so in particular when summed over the subset of $\coll_{\e m}$ of fields containing a subfield of degree $\e$.
Recall the definition (\ref{Schanuelconst}) of $S_K(n)$. Using Siegel-Brauer's Theorem,
$\delta(K)\gg_{[K:\IQ]} |\Delta_K|^{1/(\e m)}$ from Lemma 4.5 in \cite{art2} and Lemma \ref{dyadicsum} we see that also the main term converge when summed over the subset of $\coll_{\e m}$ of fields containing a subfield of degree $\e$.\\
In the proof of Theorem \ref{th2} (but now with $\e$ replaced by $\e m$ and $\coll_\e$ replaced by the subset of $\coll_{\e m}$ consisting of fields that contain a subfield of degree $\e$) we have seen that the remaining error terms coming from (\ref{Zer}), namely (\ref{cha4red2}), (\ref{2cha4red2}) and (\ref{3cha4red2}), are covered by the error term in Theorem \ref{th5}. This completes the proof of Theorem \ref{th5}. \\
As a final remark we point out that the condition $n>\max\{6\e m-6,10\}$ could be slightly relaxed since we are summing over a thinner set than $\coll_{\e m}$.

\subsection{Proof of Theorem \ref{th3}}
Let $\beta$ be as in Theorem \ref{th3} and let $f$ be the monic minimal polynomial of $\beta$ over $K$. Thus $\deg f=n$, $f$ is irreducible over $K$ and so $f$ has exactly $n$ pairwise distinct zeros. Moreover (\ref{theorem3cha4cond2}) is equivalent to $\IQ(P_f)=K$. We have seen in (\ref{M0H}) that $M_0(f)=H(\beta)^n$. Thus as $X>0$ tends to infinity the number of elements $\beta$ counted in 
Theorem \ref{th3} is given by
\begin{alignat}1
\label{Zer3}
n|\grmnnew|+O(|\grmnmnew1|)+O(|\Mred(n,X^n)|).
\end{alignat}
From (\ref{th5eq1}) and (\ref{th5eq2}), but now with $K$ of degree $\e$ instead of $\e m$, we get as $X>0$ tends to infinity
\begin{alignat*}1
|\grmnnew|&=V_{\IR}(n)^{r_K}V_{\IC}(n)^{s_K}S_K(n)X^{\e n (n+1)}\\
&+O(\delta(K)^{-\frac{\e}{2}(n-2\e+3)+1.05}X^{\e n(n+1)-n}).
\end{alignat*}
The error term above is not larger than the error term in Theorem \ref{th3}.
For the first error term in (\ref{Zer3}) we refer to (\ref{ineq620}) and then we use (\ref{heightcomp1}). In this way we see that
the first error term in (\ref{Zer3}) is also covered by the error term in Theorem \ref{th3}.
Finally due to (\ref{FS1}) the last error term in (\ref{Zer3}) is also covered by the error term in Theorem \ref{th3}.
This completes the proof of Theorem \ref{th3}.

\bibliographystyle{amsplain}
\bibliography{literature}

\end{document}